\documentclass[a4paper,11pt]{article}

\usepackage[utf8]{inputenc}
\usepackage[english,activeacute]{babel}
\usepackage{amsmath,amsthm,amsfonts,amssymb}
\usepackage{mathpartir}
\usepackage{ stmaryrd }
\usepackage{graphics}
\usepackage{enumerate}
\usepackage{mathrsfs}
\usepackage{rotating}
\input{xy}
\usepackage[all]{xy}

\newtheorem{thm}{Theorem}[subsection]
\newtheorem{cor}[thm]{Corollary}
\newtheorem{lemma}[thm]{Lemma}
\newtheorem{proposition}[thm]{Proposition}
\newtheorem{rmk}[thm]{Remark}
\theoremstyle{plain}
\newtheorem{defs}[thm]{Definition}
\theoremstyle{remark}

\input{diagxy}
\xyoption{curve}

\newbox\anglebox % large pullback angle
\setbox\anglebox=\hbox{\xy \POS(75,0)\ar@{-} (0,0) \ar@{-} (75,75)\endxy}
 
\newbox\angleboxr % reverse large pullback angle
\setbox\angleboxr=\hbox{\xy \POS(0,0)\ar@{-} (0,75) \ar@{-} (75,0)\endxy}
 
\newbox\sanglebox % small pullback angle
\setbox\sanglebox=\hbox{\xy \POS(50,0)\ar@{-} (0,0) \ar@{-} (50,50)\endxy}
 
\newbox\sangleboxr % small reverse pullback angle
\setbox\sangleboxr=\hbox{\xy \POS(0,0)\ar@{-} (0,50) \ar@{-} (50,0)\endxy}
 
\newbox\sangleboxf % small flipped pullback angle
\setbox\sangleboxf=\hbox{\xy \POS(50,50)\ar@{-} (50,0) \ar@{-} (0,50)\endxy}
 
\newbox\angleboxf % flipped pullback angle
\setbox\angleboxf=\hbox{\xy \POS(75,75)\ar@{-} (75,0) \ar@{-} (0,75)\endxy}
 
\newbox\sangleboxfr % small flipped reverse pullback angle
\setbox\sangleboxfr=\hbox{\xy \POS(0,50)\ar@{-} (50,50) \ar@{-} (0,0)\endxy}
 
\newbox\angleboxfr % small flipped reverse pullback angle
\setbox\angleboxfr=\hbox{\xy \POS(0,75)\ar@{-} (75,75) \ar@{-} (0,0)\endxy}

\title{A complete axiomatization of infinitary first-order intuitionistic logic over $\mathcal{L}_{\kappa^+, \kappa}$}
\author{Christian Esp\'indola}

\begin{document}
\date{}
\maketitle

\begin{abstract}
Given a weakly compact cardinal $\kappa$, we give an axiomatization of intuitionistic first-order logic over $\mathcal{L}_{\kappa^+, \kappa}$ and prove it is sound and complete with respect to Kripke models. As a consequence we get the disjunction and existence properties for that logic. This generalizes the work of Nadel in \cite{nadel} for intuitionistic logic over $\mathcal{L}_{\omega_1, \omega}$. When $\kappa$ is a regular cardinal such that $\kappa^{<\kappa}=\kappa$, we deduce, by an easy modification of the proof, a complete axiomatization of intuitionistic first-order logic over $\mathcal{L}_{\kappa^+, \kappa, \kappa}$, the language with disjunctions of at most $\kappa$ formulas, conjunctions of less than $\kappa$ formulas and quantification on less than $\kappa$ many variables. In particular, this applies to any regular cardinal under the Generalized Continuum Hypothesis.
\end{abstract}

\noindent $\mathbf{Keywords:}$ infinitary logics, completeness theorems, intuitionistic logic

\section{Introduction}

Completeness theorems for infinitary intuitionistic logics begun to be studied by the end of the 1970's. While systems for fragments of full first-order infinitary intuitionistic logic were considered by Makkai in \cite{makkai} and the author in \cite{espindola}, the completeness for the full first-order case has been first considered by Nadel, who in \cite{nadel} developed infinitary intuitionistic logic for countable many conjunctions/disjunctions and finite quantification, proving completeness with respect to the infinitary version of Kripke semantics. His results essentially show the completeness of intuitionistic $\mathcal{L}_{\omega_1, \omega}$ when augmented with the following axiom:

$$\bigwedge_{i<\gamma} (\phi \vee \psi_i) \to \phi \vee \bigwedge_{i<\gamma} \psi_i \qquad (1)$$
\noindent

Axiom $(1)$ is forced in every Kripke model, but is in general not valid in Heyting algebras; hence it has to be adopted. The purpose of this paper is to provide a complete axiomatization of intuitionistic logic over the language $\mathcal{L}_{\kappa^+, \kappa}$ whenever $\kappa$ is a weakly compact cardinal, with respect to Kripke models. We will see that the result of Nadel is not particular of $\mathcal{L}_{\omega_1, \omega}$ but in fact holds for any language $\mathcal{L}_{\kappa^+, \kappa}$ where $\kappa$ is weakly compact, provided the system is suitably augmented with intuitionistic rules valid in all Kripke models, which for $\kappa=\omega$ turn out to be provable from the other axioms. Hence, these results are a direct generalization of the results in \cite{nadel}. As a consequence we obtain the disjunction and existence properties for our logic, supporting its constructive character. Also, by restricting the language to $\mathcal{L}_{\kappa^+, \kappa, \kappa}$ (that is, only allowing conjunctions of less than $\kappa$ many formulas), the resulting axiomatization reduces to one that is complete, for Kripke semantics, whenever $\kappa^{<\kappa}=\kappa$, not necessarily weakly compact.

Finally, we should mention that the author has obtained in \cite{espindola} a strong completeness theorem for intuitionistic logic over $\mathcal{L}_{\kappa, \kappa}$ whenever $\kappa$ is a weakly compact cardinal. A pertinent observation is that in this paper whenever we talk about completeness we always refer to weak completeness, i.e., valid sentences being provable (in the empty theory). This is of course equivalent to having completeness for theories axiomatized by at most $\kappa$ many axioms, since the conjunction of those axioms is expressible in the language of $\mathcal{L}_{\kappa^+, \kappa}$. As noted already in \cite{espindola}, this is best possible in terms of the cardinality of the theories involved, since we will see that there are theories with $\kappa^+$ many axioms that are consistent but have no Kripke model.

This work is a continuation of the investigation begun in \cite{espindola} on infinitary first-order categorical logic, but the proofs will only use lattice-theoretic methods on this ocassion. The naturality of the given axiomatization, of course, is motivated by the category-theoretic interpretation. The structure is as follows. We first introduce the system for infinitary intuitionistic logic over $\mathcal{L}_{\kappa^+, \kappa}$. We then work with a coherent fragment of our system and prove it is complete with respect to set-valued models whenever $\kappa$ is weakly compact. This is based on a development of an infinitary version of Priestley duality for certain lattices, together with a topological approach to a Rasiowa-Sikorski type result. As a consequence, we obtain also the completeness of the full system with respect to validity in Kripke models.

\subsection{Infinitary first-order systems}

Let $\kappa$ be a regular cardinal such that $\kappa^{<\kappa}=\kappa$. The syntax of $\kappa$-first-order logic consists of a (well-ordered) set of sorts and a set of function and relation symbols, these latter together with the corresponding type, which is a subset with less than $\kappa$ many sorts. Therefore, we assume that our signature may contain relation and function symbols on $\gamma<\kappa$ many variables, and we suppose there is a supply of $\kappa$ many fresh variables of each sort. 
 
We use sequent style calculus to formulate the axioms of $\kappa$-first-order logic, as can be found, e.g., in \cite{johnstone}, D1.3. A sequent $\phi \vdash_{\mathbf{x}} \psi$ is however equivalent to the formula $\forall \mathbf{x} (\phi \to \psi)$, though we find it convenient to keep the sequent notation since it is closer to the categorical semantics that will be used in the proofs. The system is described in the following:

\begin{defs}\label{sfol}
 The system of axioms and rules for $\kappa$-first-order logic consists of

\begin{enumerate}
 \item Structural rules:
 \begin{enumerate}
 \item Identity axiom:
\begin{mathpar}
\phi \vdash_{\mathbf{x}} \phi 
\end{mathpar}
\item Substitution rule:
\begin{mathpar}
\inferrule{\phi \vdash_{\mathbf{x}} \psi}{\phi[\mathbf{s}/\mathbf{x}] \vdash_{\mathbf{y}} \psi[\mathbf{s}/\mathbf{x}]} 
\end{mathpar}

where $\mathbf{y}$ is a string of variables including all variables occurring in the string of terms $\mathbf{s}$.
\item Cut rule:
\begin{mathpar}
\inferrule{\phi \vdash_{\mathbf{x}} \psi \\ \psi \vdash_{\mathbf{x}} \theta}{\phi \vdash_{\mathbf{x}} \theta} 
\end{mathpar}
\end{enumerate}

\item Equality axioms:

\begin{enumerate}
\item 

\begin{mathpar}
\top \vdash_{x} x=x 
\end{mathpar}

\item 

\begin{mathpar}
(\mathbf{x}=\mathbf{y}) \wedge \phi[\mathbf{x}/\mathbf{w}] \vdash_{\mathbf{z}} \phi[\mathbf{y}/\mathbf{w}]
\end{mathpar}

where $\mathbf{x}$, $\mathbf{y}$ are contexts of the same length and type and $\mathbf{z}$ is any context containing $\mathbf{x}$, $\mathbf{y}$ and the free variables of $\phi$.
\end{enumerate}

\item Conjunction axioms and rules:

$$\bigwedge_{i<\gamma} \phi_i \vdash_{\mathbf{x}} \phi_j$$

\begin{mathpar}
\inferrule{\{\phi \vdash_{\mathbf{x}} \psi_i\}_{i<\gamma}}{\phi \vdash_{\mathbf{x}} \bigwedge_{i<\gamma} \psi_i}
\end{mathpar}

for each cardinal $\gamma<\kappa^+$.

\item Disjunction axioms and rules:

$$\phi_j \vdash_{\mathbf{x}} \bigvee_{i<\gamma} \phi_i$$

\begin{mathpar}
\inferrule{\{\phi_i \vdash_{\mathbf{x}} \theta\}_{i<\gamma}}{\bigvee_{i<\gamma} \phi_i \vdash_{\mathbf{x}} \theta}
\end{mathpar}

for each cardinal $\gamma<\kappa^+$.

\item Implication rule:
\begin{mathpar}
\mprset{fraction={===}}
\inferrule{\phi \wedge \psi \vdash_{\mathbf{x}} \eta}{\phi \vdash_{\mathbf{x}} \psi \to \eta}
\end{mathpar}

\item Existential rule:
\begin{mathpar}
\mprset{fraction={===}}
\inferrule{\phi \vdash_{\mathbf{x} \mathbf{y}} \psi}{\exists \mathbf{y}\phi \vdash_{\mathbf{x}} \psi}
\end{mathpar}

where no variable in $\mathbf{y}$ is free in $\psi$.

\item Universal rule:
\begin{mathpar}
\mprset{fraction={===}}
\inferrule{\phi \vdash_{\mathbf{x} \mathbf{y}} \psi}{\phi \vdash_{\mathbf{x}} \forall \mathbf{y} \psi}
\end{mathpar}

where no variable in $\mathbf{y}$ is free in $\phi$.

\item Small distributivity axiom:

$$\bigwedge_{i<\gamma} \phi \vee \psi_i \vdash_{\mathbf{x}} \phi \vee \bigwedge_{i<\gamma} \psi_i$$

for each cardinal $\gamma<\kappa^+$.

\item Dual distributivity rule:

\begin{mathpar}
\inferrule{\bigwedge_{g \in \gamma^{\beta+1}, g|_{\beta}=f} \phi_{g} \vdash_{\mathbf{x}} \phi_{f} \\ \beta<\kappa, f \in \gamma^{\beta} \\\\ \phi_{f} \dashv \vdash_{\mathbf{x}} \bigvee_{\alpha<\beta}\phi_{f|_{\alpha}} \\ \beta < \kappa, \text{ limit }\beta, f \in \gamma^{\beta}}{\bigwedge_{f \in B} \bigvee_{\beta<\delta_f}\phi_{f|_{\beta+1}} \vdash_{\mathbf{x}} \phi_{\emptyset}}
\end{mathpar}
\\
for each cardinal $\gamma<\kappa^+$. Here $B \subseteq \gamma^{< \kappa}$ consists of the minimal elements of a given bar\footnote{A bar over the tree $\gamma^{\kappa}$ is an upward closed subset of nodes intersecting every branch of the tree.} over the tree $\gamma^{\kappa}$, and the $\delta_f$ are the levels of the corresponding $f \in B$.

\item Transfinite transitivity rule:

\begin{mathpar}
\inferrule{\phi_{f} \vdash_{\mathbf{y}_{f}} \bigvee_{g \in \gamma^{\beta+1}, g|_{\beta}=f} \exists \mathbf{x}_{g} \phi_{g} \\ \beta<\kappa, f \in \gamma^{\beta} \\\\ \phi_{f} \dashv \vdash_{\mathbf{y}_{f}} \bigwedge_{\alpha<\beta}\phi_{f|_{\alpha}} \\ \beta < \kappa, \text{ limit }\beta, f \in \gamma^{\beta}}{\phi_{\emptyset} \vdash_{\mathbf{y}_{\emptyset}} \bigvee_{f \in B}  \exists_{\beta<\delta_f}\mathbf{x}_{f|_{\beta +1}} \bigwedge_{\beta<\delta_f}\phi_{f|_{\beta+1}}}
\end{mathpar}
\\
for each cardinal $\gamma<\kappa^+$, where $\mathbf{y}_{f}$ is the canonical context of $\phi_{f}$, provided that, for every $f \in \gamma^{\beta+1}$,  $FV(\phi_{f}) = FV(\phi_{f|_{\beta}}) \cup \mathbf{x}_{f}$ and $\mathbf{x}_{f|_{\beta +1}} \cap FV(\phi_{f|_{\beta}})= \emptyset$ for any $\beta<\gamma$, as well as $FV(\phi_{f}) = \bigcup_{\alpha<\beta} FV(\phi_{f|_{\alpha}})$ for limit $\beta$. Here $B \subseteq \gamma^{< \kappa}$ consists of the minimal elements of a given bar over the tree $\gamma^{\kappa}$, and the $\delta_f$ are the levels of the corresponding $f \in B$.

\end{enumerate}
\end{defs}
\noindent
The motivation of the transfinite transitivity rule has been explained at length in \cite{espindola}, to which we refer the reader.

\section{Completeness of intuitionistic logic over $\mathcal{L}_{\kappa^+, \kappa}$ for weakly compact $\kappa$}

Throughout this section we will always assume that $\kappa$ is a weakly compact cardinal. The $\kappa$-coherent fragment of $\kappa$-first-order logic over $\mathcal{L}_{\kappa^+, \kappa}$ consists of sequents of coherent formulas, which are inductively built using conjunctions and disjunctions of at most $\kappa$-many subformulas, and existential quantification on less than $\kappa$ many formulas. The axioms and rules are those involving only these quantifiers and connectives, and instantiated only on $\kappa$-coherent formulas.

\subsection{$\kappa$-distributive lattices}

In \cite{espindola} $\kappa$-distributive lattices were defined as an infinitary generalization of distributive lattices in such a way that they have the right property to set up a Priestley duality theory for them. More especifically, the distributivity property was chosen so as to allow the existence of enough $\kappa$-complete, $\kappa$-prime filters\footnote{In a $\kappa$-complete lattice, a filter is $\kappa$-complete if it is closed under intersection of less than $\kappa$ many elements, and it is $\kappa$-prime if whenever a join of less than $\kappa$ many elements is in the filter then one of them must be in it.}, in the sense that they separate points (that is, whenever $a \nleq b$, there is a $\kappa$-complete, $\kappa$-prime filter containing $a$ but not $b$).

To motivate the definition, suppose we have a $\kappa$-complete lattice with the separation property just mentioned. Consider, for $\gamma<\kappa$, a tree $\gamma^{<\gamma}$ (a $\gamma$-branching tree of height $\gamma$) whose nodes are elements of the lattice, and with the opposite order. Assume that the tree satisfies the following two conditions:

\begin{itemize}
\item the element in every node of the tree is below the join of the elements at its immediate successors

\item the element in a node at a limit level is the meet of its predecessors

\end{itemize}

If the $\kappa$-complete, $\kappa$-prime filters separate points, then the lattice must satisfy also the following property: the element at the root of the tree is below the join of the minimal elements at any given level of the tree (not necessarily finite). Indeed, if that was not the case, there would be a $\kappa$-complete, $\kappa$-prime filter $P$ sending the element at the root to $1$ and the join of the elements at some given level $\beta$ of the tree to $0$. But then we could inductively define a cofinal branch in the tree composed of elements that are sent to $1$ by $P$, and so the element in this branch at level $\beta$ would have to be one of these, contradicting the choice of $P$. This motivates the following:

\begin{defs}
A $\kappa$-complete lattice with possibly a set of joins of $\kappa$ many elements will be called $\kappa$-distributive if for each $\gamma<\kappa$ and all elements $\{a_f: f \in \gamma^{\beta}, \beta<\kappa\}$  such that 

$$a_{f} = \bigvee_{g \in \gamma^{\beta+1}, g|_{\beta}=f} a_{g}$$
\noindent for all $f \in \gamma^{\beta}, \beta<\kappa$, and 

$$a_{f} = \bigwedge_{\alpha<\beta}a_{f|_{\alpha}}$$
\noindent for all limit $\beta$, $f \in \gamma^{\beta}, \beta<\kappa$, we have that the join $\bigvee_{f \in B} \bigwedge_{\beta<\delta_f}a_{f|_{\beta}}$ exists and is equal to $a_{\emptyset}$, where $B \subseteq \gamma^{< \kappa}$ consists of the minimal elements of any given bar over the tree $\gamma^{<\kappa}$.

The lattice is called $\kappa^+$-distributive if we can take $\gamma=\kappa$ in the definition above.

\end{defs}

This strong distributivity property is, as we will see, enough to guarantee the separation property of the lattice even if we restrict ourselves to those $\kappa$-complete, $\kappa$-prime filters that also preserve joins of $\kappa$ many elements: 

\begin{proposition}\label{filter}
Let $\kappa$ be a regular cardinal such that $\kappa^{<\kappa}=\kappa$. Then, any non-trivial $\kappa$-complete, $\kappa^+$-distributive lattice of cardinality at most $\kappa$, with a set $\mathcal{S}$ of at most $\kappa$ joins of $\kappa$ many elements each, has a $\kappa$-complete, $\kappa$-prime filter preserving also any join in $\mathcal{S}$. Moreover, given $a \nleq b$, there exists a $\kappa$-complete, $\kappa$-prime filter preserving also any join in $\mathcal{S}$, and containing $a$ but not $b$ (i.e., $\kappa$-complete, $\kappa$-prime filters that preserve also the joins in $\mathcal{S}$ separate points).
\end{proposition}

\begin{proof}
Let $\mathcal{L}$ be a $\kappa$-complete, $\kappa^+$-distributive lattice of cardinality at most $\kappa$. We will first show how to build a $\kappa$-complete, $\kappa$-prime filter $\mathcal{F}$ in $\mathcal{L}$ preserving also any join of $\kappa$ many elements from a set of at most $\kappa$ many of these.

Start with a well-ordering $f: \kappa \times \kappa \to \kappa$ with the property that $f(\beta, \gamma) \geq \gamma$. For example, the canonical well-ordering defined by induction on $\max (\beta, \gamma)$ as follows:
 
 $$f(\beta, \gamma)=\begin{cases} \sup\{f(\beta', \gamma')+1: \beta', \gamma'<\gamma\}+\beta  & \mbox{if }  \beta<\gamma  \\ \sup\{f(\beta', \gamma')+1: \beta', \gamma'<\beta\}+\beta+\gamma & \mbox{if } \gamma \leq \beta \end{cases}$$
 
\noindent clearly satisfies the required property.
 
 For each $a \in \mathcal{L}$, let $\mathcal{C}(a)$ be the set of all tuples $(b_\alpha)_{\alpha<\lambda}$ of elements of $\mathcal{L}$ such that:

\begin{itemize}

\item $\lambda \leq \kappa$ and if $\lambda=\kappa$ then $\bigvee_{\alpha<\kappa} b_\alpha$ is in $\mathcal{S}$

\item $a = \bigvee_{\alpha<\lambda}b_{\alpha}$.

\end{itemize}
Assume without loss of generality that $\mathcal{C}(a)$ is well-ordered and has order type $\kappa$, repeating tuples, if needed (since $\kappa^{<\kappa}=\kappa$, this is clearly possible). Then we can build a tree of height $\kappa$ whose nodes are elements of $\mathcal{L}$ through the following transfinite recursive procedure. The root of the tree is the top element of $\mathcal{L}$. Assuming that the tree is defined for all levels $\lambda<\mu$; we show how to define the nodes of level $\mu$. Suppose first that $\mu$ is a successor ordinal $\mu=\alpha+1$, and let $\alpha=f(\beta, \gamma)$. Since by hypothesis $f(\beta, \gamma) \geq \gamma$, the nodes $\{p_i\}_{i<m_{\gamma}}$ at level $\gamma$ are defined. To define the nodes at level $\alpha+1$, we need to define the successors of a node $n$ there; for this purpose, take then the $\beta-th$ tuple $(b_\alpha)_{\alpha<\lambda} \in \mathcal{C}(p)$ over the predecessor $p$ of $n$ at level $\gamma$, and define the successors of $n$ to be $(n \wedge b_\alpha)_{\alpha<\lambda}$. Suppose now that $\mu$ is a limit ordinal. Then define every node at level $\mu$ to be the conjunction of the predecessors.
 
 By construction, and because of the distributivity property, the join of all elements along the minimal elements of any given bar over the subtree over any given node $c$ is equal to $c$, and hence, the join of all elements along the minimal elements of any given bar over the tree is $1$. Therefore, at least one node in that join is not $0$. This readily implies that there must exist a cofinal branch $\mathbb{B}$ in the tree entirely composed of nodes that are not $0$. Then we can define $\mathcal{F}$ by stipulating $a \in \mathcal{F}$ if and only if there is some node $b$ in $\mathbb{B}$ such that $b \leq a$. 
 
 $Claim:$ $\mathcal{F}$ is a $\kappa$-complete, $\kappa$-prime filter preserving all the joins in $\mathcal{S}$.
 
 \begin{enumerate}
  \item It is clearly seen to contain $1$ and not to contain $0$.
  \item If $a \leq b$ and $a \in \mathcal{F}$ then clearly $b \in \mathcal{F}$.
  \item It is closed under $\gamma<\kappa$ conjunctions. Indeed, if for $\alpha<\nu$ we have that $b_{\alpha}$ at level $l(\alpha)$ witnesses that $a_{\alpha} \in \mathcal{F}$, then $\bigwedge_{\alpha<\nu}b_{\alpha}$ will be the node of the branch at level $\sup_{\alpha<\nu}l(\alpha)$, and will witness that $\bigwedge_{\alpha<\nu}a_{\alpha} \in \mathcal{F}$.
  \item For the primeness property, suppose that $b$ witnesses that $\bigvee_{\alpha<\gamma}a_{\alpha}$ is in $\mathcal{F}$. Then we have $(b \wedge a_{\alpha})_{\alpha<\gamma} \in \mathcal{C}(b)$, so that by construction, for some level in the subtree over the node $b$, every element $c$ in that level will have successors of the form $b \wedge a_{\alpha} \wedge c$, which implies, by definition, that some $a_{\alpha}$ will be in $\mathcal{F}$.  
 \end{enumerate}

Let us now prove that given $a \nleq b$, we can construct $\mathcal{F}$ in such a way that it contains $a$ but not $b$. First, construct by transfinite recursion a tree exactly as before but using $a$ as a root. Since $a \nleq b$, this means, with a similar argument to the one we gave before, that there is a cofinal branch $\mathbb{B}$ composed of nodes $c \nleq b$. We can then build $\mathcal{F}$ as before, and since $a \in \mathbb{B}$, $a$ will be in $\mathcal{F}$, but $b$ will not. This finishes the proof.

\end{proof}

The following corollary and (slightly rewritten) proof had been obtained already in \cite{espindola}:

\begin{cor}\label{filterwc}
Let $\kappa$ be a weakly compact cardinal. Then, in any non-trivial $\kappa$-complete, $\kappa$-distributive lattice of cardinality at most $\kappa$, for any $\kappa$-complete filter $\mathcal{F}$ disjoint from a $\kappa$-complete ideal $\mathcal{I}$, there exists a $\kappa$-complete, $\kappa$-prime filter containing $\mathcal{F}$ and disjoint from $\mathcal{I}$. 
\end{cor}

\begin{proof}
Define a new propositional variable $P_a$ for each element $a$ of $\mathcal{L}$ and consider the theory axiomatized by the following axioms: $\top \vdash P_a$ (for all $a \in \mathcal{F}$); $P_a \vdash \bot$ (for all $a \in \mathcal{I}$); $P_a \vdash P_b$ (for each pair $a \leq b$); $\bigwedge_{i \in I}P_{a_i} \vdash P_{\left(\bigwedge_{i \in I}a_i\right)}$ (for all $a_i$ and $|I|<\kappa$); and $P_{\left(\bigvee_{i \in I}a_i\right)} \vdash \bigvee_{i \in I}P_{a_i}$ (for all $a_i$ and $|I|<\kappa$). This is a theory over the $\kappa$-coherent fragment which has cardinality at most $\kappa$. Each subtheory of cardinality less than $\kappa$ involves $\gamma<\kappa$ propositional variables, whose correspondent elements generate a (non trivial) $\kappa$-complete sublattice $\mathcal{L}'$ of $\mathcal{L}$ of cardinality $2^{\gamma}<\kappa$ (as every element there is equivalent to one of the form $\bigvee_{i<\gamma} \bigwedge_{j<\gamma} P_{a_{ij}}$). If $a$ is the intersection of all elements in $\mathcal{F} \cap \mathcal{L}'$ and $b$ the union of all elements in $\mathcal{I} \cap \mathcal{L}'$, by an entirely analogous proof to that of Proposition \ref{filter}, there is a $\kappa$-complete, $\kappa$-prime filter in $\mathcal{L}'$ containing $a$ but not $b$. This shows that each subtheory of cardinality less than $\kappa$ has a model. Since $\kappa$ is weakly compact, the whole theory has a model, which corresponds to a $\kappa$-complete, $\kappa$-prime filter of $\mathcal{L}$ containing $\mathcal{F}$ and disjoint from $\mathcal{I}$.
\end{proof} 

Corollary \ref{filterwc} allows to generalize the theory of Priestley duality for distributive lattices to the infinitary case. More precisely, define a $\kappa$-space as a space with a $\kappa$-topology, i.e., a topology where the intersection of less than $\kappa$ many open sets is open. A subbasis of a $\kappa$-space is a collection of open sets such that any open is a union of a $\kappa$-small intersection of sets from the collection. A $\kappa$-Stone space as a Hausdorff $\kappa$-space with a basis of clopen sets that is $\kappa$-compact: every open cover has a $\kappa$-small subcover. We have the following:

\begin{lemma}
Let $\kappa$ be a weakly compact cardinal. If every subbasic cover of a $\kappa$-space $X$ with a subbasis of cardinality at most $\kappa$ has a $\kappa$-small subcover, $X$ is $\kappa$-compact.
\end{lemma}

\begin{proof}
This is the infinitary version of the well known Alexander's subbase lemma, and the proof is essentially contained in \cite{ms}. Let $\mathcal{L}$ be the $\kappa$-complete, $\kappa$-distributive sublattice of $\mathcal{P}(X)$ generated by the elements of the subbasis, and let $\mathcal{B}$ be the $\kappa$-complete ideal $\mathcal{L}$ of generated by the elements of the basis. Assume that every subbasic cover admits a $\kappa$ subcover but there is a family of basic open with no basic $\kappa$-subcover that covers $X$. We will show that the family does not cover $X$ either. By hypothesis, $X$ is not in $\mathcal{B}$, and hence there is a proper $\kappa$-complete, $\kappa$-prime ideal $\mathcal{I}$ containing $\mathcal{B}$.  Then no $\kappa$-subfamily of subbasic opens in $\mathcal{I}$ covers $X$, and so by hypothesis the union of all subbasic opens in $\mathcal{I}$ does not cover $X$ either. Let $x \in X$ not in this union. We claim that $x$ is not in the union of the basic opens either. Indeed, otherwise $x$ would belong to an intersection of subbasic opens, and since this latter is in $\mathcal{I}$, which is $\kappa$-complete, at least one of the subbasic opens of the intersection would be in $\mathcal{I}$, contrary to the choice of $x$. 
\end{proof}

\begin{lemma}
The space of $\kappa$-complete, $\kappa$-prime filters of a $\kappa$-distributive lattice is a $\kappa$-Stone space. 
\end{lemma}

\begin{proof}
Straightforward generalization of the finitary case, making use of Proposition \ref{filterwc} and the infinitary version of Alexander's subbase lemma.
\end{proof}

A $\kappa$-Priestley space is a $\kappa$-Stone space with a partial order $\leq$ satisfying the usual separation axiom: if $a \nleq b$, there is a clopen upset containing $a$ but not $b$. Using once more Proposition \ref{filter}, we get:

\begin{lemma}
The space of $\kappa$-complete, $\kappa$-prime filters of a $\kappa$-distributive lattice with the partial order given by inclusion is a $\kappa$-Priestley space. 
\end{lemma}

\begin{thm}\label{priestley}
There is a duality between the category of $\kappa$-distributive lattices and the category of $\kappa$-Priestley spaces, given by homming into $2$.
\end{thm}

Priestley duality is useful, for example, in providing a topological proof of Rasiowa-Sikorski lemma for distributive lattices. This was done in \cite{goldblatt}, and it turns out that this proof is adequate for its generalization to the infinitary case. Before getting to it, we need the following infinitary version of Baire category theorem for $\kappa$-spaces: the intersection of $\kappa$ many dense open sets is dense. The proof of this latter statement, in turn, is essentially contained in \cite{sikorski}:

\begin{proposition}\label{baire}
\textbf{($\kappa$-Baire category theorem)} Let $\kappa$ be a weakly compact cardinal. Then in the $\kappa$-space of $\kappa$-complete, $\kappa$-prime filters of a $\kappa$-complete, $\kappa$-distributive lattice of cardinality at most $\kappa$, let a family of $\kappa$ many dense open sets be given, in such a way that the intersection of less than $\kappa$ many of the sets in the family is nonempty. Then the intersection of all the sets in the family is dense.
\end{proposition}

\begin{proof}
The $\kappa$ space $\mathcal{P}$ of $\kappa$-complete, $\kappa$-prime filters has a basis of cardinality $\kappa$ consisting of clopen sets of the form $\phi(a) \cap \phi(b)^c$, where $\phi(z)$ is the set of $\kappa$-complete, $\kappa$-prime filters containing the point $z$. By theorem $(x)$ of \cite{sikorski}, $\mathcal{P}$ is homeomorphic to a subset of the generalized Cantor space $2^{\kappa}$. Moreover, since $\mathcal{P}$ is $\kappa$-compact and $2^{\kappa}$ is Hausdorff (in fact, it is a complete $\kappa$-metric space, cf. \cite{sikorski}), such a subset is closed, and therefore complete. Since the intersection of less than $\kappa$ many of the dense open sets is nonempty, the proof of theorem $(xv)$ there applies mutatis mutandis to this subspace (the proviso is necessary to be able to continue the definition of the $\kappa$-sequence beyond limit ordinals), showing that $\mathcal{P}$ satisfies the $\kappa$-Baire category theorem.
\end{proof}

We now turn to the Rasiowa-Sikorski type result. In a $\kappa$-distributive lattice $\mathcal{L}$, the join $\bigvee_{i<\kappa}a_i$ is called \emph{distributive} if for any $c \in \mathcal{L}$ we have $c \wedge \bigvee_{i<\kappa}a_i = \bigvee_{i<\kappa} (c \wedge b_i)$. Dually, the meet $\bigwedge_{i<\kappa}b_i$ is called distributive if for any $c \in \mathcal{L}$ we have $\bigwedge_{i<\kappa} (c \vee b_i) = c \vee \bigwedge_{i<\kappa}b_i$. We can now deduce the following:

\begin{proposition}\label{rs}
Let $\kappa$ be a weakly compact cardinal. Let $\mathcal{L}$ be a $\kappa$-distributive lattice and suppose that $\{c_i: i<\kappa\}$ (respectively, $\{d_i: i<\kappa\}$) are sets of joins (respectively, meets) of $\kappa$-many elements ($c_i=\bigvee_{j<\kappa}a_{ij}$ for each $i<\kappa$, $d_i=\bigwedge_{j<\kappa}b_{ij}$ for each $i<\kappa$) that are distributive. If for each $\gamma<\kappa$ there is a $\kappa$-complete, $\kappa$-prime filter preserving the joins and meets $\{c_i: i<\gamma\}$ and $\{d_i: i<\gamma\}$, then, given $a \nleq b$, there is a $\kappa$-complete, $\kappa$-prime filter preserving all the joins and meets and containing $a$ but not $b$.
\end{proposition}

\begin{proof}
The proof is a straightforward generalization of the proof for the finitary case given in \cite{goldblatt}, which made use of Priestley duality and the Baire category theorem. In our case, we can use the duality from Theorem \ref{priestley} and Proposition \ref{baire}. It is easy to check, using distributivity, that a similar proof to the one in \cite{goldblatt} shows here that the subset of $\kappa$-complete, $\kappa$-prime filters preserving a given join or a given meet is open and dense. Since the subset of $\kappa$-complete, $\kappa$-prime filters preserving the joins and meets $\{c_i: i<\gamma\}$ and $\{d_i: i<\gamma\}$, for $\gamma<\kappa$, is nonempty, applying Proposition \ref{baire} we get the result.
\end{proof}

\subsection{Completeness of first-order intuitionistic logic over $\mathcal{L}_{\kappa^+, \kappa}$ and of the $\kappa^+$-coherent fragment}

We are now going to apply to apply Proposition \ref{rs} to the lattice of formulas (not necessarily closed) of the $\kappa^+$-coherent fragment of logic. This fragment consists of sequents of $\kappa^+$-coherent formulas, and these latter are formed using only conjunctions and disjunctions of size at most $\kappa$ (provided the set of free variables of te resulting formulas has size less than $\kappa$) and existential quantification over less than $\kappa$ many variables. The axioms are those of the full first-order fragment that only involve these connectives and quantifiers, and where the instantiations in schemata is done on $\kappa^+$-coherent formulas only.

Consider the full Lindenbaum-Tarski algebra of all (not necessarily closed) $\kappa^+$-coherent formulas, where the order relation is given by $\phi(\mathbf{x}) \leq \psi(\mathbf{y})$ if and only if the sequent $\phi(\mathbf{x}) \vdash_{\mathbf{x}\mathbf{y}} \psi(\mathbf{y})$ is provable in a previously given theory of cardinality less than $\kappa$. This is a $\kappa$-complete lattice. Let $\phi \vdash_{\mathbf{x}} \psi$ be a sequent which is not provable. Consider now the set of at most $\kappa$ joins (respectively, meets) $\{c_i=\bigvee_{j<\kappa}a_{ij}: i<\kappa\}$ (respectively, $\{d_i=\bigwedge_{j<\kappa}b_{ij}: i<\kappa\}$) that appear as subformulas of the axioms of the theory or of $\phi$ or $\psi$, considering for each existential subformula $\exists \mathbf{x} \phi(\mathbf{x})$, also the corresponding join over all formulas $\phi(\mathbf{y})$ when $\mathbf{y}$ ranges over all tuples of variables. Let $\mathcal{L}$ be the $\kappa$-complete sublattice of cardinality at most $\kappa$ generated by the subformulas of the axioms of the theory or of $\phi$ or $\psi$, as well as the substitution of these for terms. It is clear that the set of joins and meets considered is distributive.

  Clearly, $\mathcal{L}$ is $\kappa$-distributive, but we will prove even that it is $\kappa^+$-distributive, or equivalently, that the $\kappa$-complete, $\kappa$-prime filters preserving the given set of $\kappa^+$-joins separate $\phi$ and $\psi$ (since these latter are arbitrary). For this purpose it is in turn enough to verify the following condition used in the proof of Proposition \ref{filter}: if $\phi \nleq \psi$, the tree built over $\phi$ by transfinite recursion has a cofinal branch composed of formulas $\theta \nleq \psi$. This condition will be a consequence of the transfinite transitivity rule. Indeed, it suffices to verify that for every bar $b$, if each formula $\bigwedge_{i<\delta_f} \phi_{f|_i}$ entails $\psi$, so does $\phi$. But if we consider $\phi_f=\exists \mathbf{x} \theta(\mathbf{x})$, without loss of generality, as a $\kappa^+$-join of the formulas of the form $\phi_g=\theta(\mathbf{y})$ when $\mathbf{y}$ ranges over all tuples of variables not appearing free in $\psi$ or in $\mathbf{x}$, then we will have $\phi_f=\bigvee_{\mathbf{y} \cap FV(\psi)=\emptyset}\exists \mathbf{y} \theta(\mathbf{y})$, and by transfinite transitivity we get that $\phi$ entails the join of the formulas $\exists \mathbf{y}\mathbf{y'}...\bigwedge_{i<\delta_f} \phi_{f|_i}$. It is now immediate that each $\exists \mathbf{y}\mathbf{y'}...\bigwedge_{i<\delta_f} \phi_{f|_i}$ entails $\psi$, since the variables $\mathbf{y}\mathbf{y'}...$ do not appear free in $\psi$, and so they can be moved to the context by the existential quantification rule. It follows that $\mathcal{L}$ is indeed $\kappa^+$-distributive.

  In order to apply Proposition \ref{rs} to $\mathcal{L}$, we need to verify the following condition: for each $\gamma<\kappa$ there is a $\kappa$-complete, $\kappa$-prime filter preserving the joins and meets $\{c_i: i<\gamma\}$ and $\{d_i: i<\gamma\}$. This condition, on the other hand, will be a consequence of the dual distributivity rule, as seen in the following:

\begin{lemma}\label{al}
Let $\kappa$ be a weakly compact cardinal. Given $a  \nleq b$ in $\mathcal{L}$, for each $\gamma<\kappa$, there is a $\kappa$-complete, $\kappa$-prime filter preserving the joins and meets $\{c_i: i<\gamma\}$ and $\{d_i: i<\gamma\}$, and containing $a$ but not $b$.
\end{lemma}

\begin{proof}
Let $\mathbf{x_i}$ be the free variables of $c_i$ and $d_i$, and define $\mathbf{x}=\cup_{i<\gamma}\mathbf{x_i} \cup FV(a) \cup FV(b)$. Consider the $\kappa$-complete, $\kappa^+$-distributive sublattice $\mathcal{L}'$ generated by $a$, $b$, $c_i$, $d_i$, $a_{ij}$ and $b_{ij}$ for $i<\gamma$, $j<\kappa$. Since all formulas in $\mathcal{L}'$ have free variables amongst $\mathbf{x}$, the dual distributivity rule implies that, besides $\mathcal{L}'$, also the dual lattice $(\mathcal{L}')^*$ is $\kappa^+$-distributive, and so by Proposition \ref{filter}, there exists a $\kappa$-complete, $\kappa$-prime filter $\mathcal{F}'$ of $\mathcal{L}'$ preserving the joins and meets $\{c_i: i<\gamma\}$ and $\{d_i: i<\gamma\}$, and containing $a$ but not $b$. We will extend $\mathcal{F}'$ to a $\kappa$-complete, $\kappa$-prime filter $\mathcal{F}$ of $\mathcal{L}$ that still preserves the given joins and meets and does not contain $b$. Let $\{d_{i_l}: l<\delta \leq \gamma\}$ be the set of those meets that do not belong to $\mathcal{F}'$. Since this latter preserves them, there are nonempty subsets $S_{i_l} \subseteq \{b_{i_lj}: j<\kappa\}$ not in $\mathcal{F}'$. Let $\{c_{i_j}: j<\delta'\}$ be the set of joins among the $\{c_i: i<\gamma\}$ which do not belong to $\mathcal{F}'$. Let $\mathcal{G}$ be the $\kappa$-complete filter of $\mathcal{L}$ generated by $\mathcal{F}'$, and let $\mathcal{I}$ be the $\kappa$-complete ideal of $\mathcal{L}$ generated by $b$, the elements in all the $S_{i_l}$ for $l<\delta$ and the $c_{i_j}$ for $j<\delta'$. Using the $\kappa$-primeness (in $\mathcal{L}'$) of $\mathcal{F}'$, it is easy to see that $\mathcal{G} \cap \mathcal{I}=\emptyset$. Since $\kappa$ is weakly compact, by Corollary \ref{filterwc} there is a $\kappa$-complete, $\kappa$-prime filter $\mathcal{F}$ of $\mathcal{L}$ containing $\mathcal{G}$ (and hence $\mathcal{F}'$) and disjoint from $\mathcal{I}$. Therefore, $\mathcal{F}$ still preserves all the joins and meets  $\{c_i: i<\gamma\}$ and $\{d_i: i<\gamma\}$.
\end{proof}

\begin{cor}\label{cohcomp}
If $\kappa$ is weakly compact, $\kappa^+$-coherent theories of cardinality at most $\kappa$ are complete with respect to set-valued models of cardinality at most $\kappa$.
\end{cor}

\begin{proof}
Let $\phi \vdash_{\mathbf{x}} \psi$ be a sequent which is not provable in the theory. Form the lattice $\mathcal{L}$ as before. Each existential formula $\exists \mathbf{x} \phi(\mathbf{x})$ is a $\kappa^+$-join of the formulas $\phi(\mathbf{t})$, where $\mathbf{t}$ ranges over all tuples of terms with free variables different from those appearing free in $\psi$. By Lemma \ref{al} and the discussion above about the applicability of Proposition \ref{rs}, there is a $\kappa$-complete, $\kappa$-prime filter $\mathcal{F}$ preserving all $\kappa^+$-joins and meets appearing in axioms of the theory or in $\phi$ or $\psi$ (including $\kappa^+$ joins arising from existential quantification), and containing $\phi$ but not $\psi$. Consider the structure whose underlying set consists of all (equivalence classes of) terms of the language, and where an atomic formula $R(\mathbf{t})$ holds if and only if $R(\mathbf{t})$ is in $\mathcal{F}$. Then such a structure is a model of the theory of cardinality at most $\kappa$ where the sequent $\phi \vdash_{\mathbf{x}} \psi$ does not hold.
\end{proof}

\begin{rmk}
Corollary \ref{cohcomp} is optimal with respect to the size of the theory. Indeed, by a result of Specker (see \cite{specker}), if $\kappa^{<\kappa}=\kappa$ there is a $\kappa^+$-Aronszajn tree, for which the theory of a cofinal branch is a consistent $\kappa^+$-coherent theory of cardinality $\kappa^+$ which has no models.
\end{rmk}

Having now at hand a completeness theorem for $\kappa^+$-coherent theories, we can adapt the proof of completeness for Kripke semantics (see \cite{dvd}). As a result, we get:

\begin{thm}\label{bt}
 If $\kappa$ is weakly compact, intuitionistic first-order theories over $\mathcal{L}_{\kappa^+, \kappa}$ of cardinality at most $\kappa$ are complete with respect to Kripke models.
\end{thm}

\begin{proof}
We generalize the usual proof of \cite{dvd} for the finitary case, which proceeds in two steps. First, given a theory $\mathbb{T}$ such that $\mathbb{T} \nvdash \phi$, there is a prime theory $\mathbb{T}'\supseteq \mathbb{T}$ in an extended language with $\mathbb{T}' \nvdash \phi$ (here a prime theory is a deductively closed theory with the disjunction and existence properties). In our case this is automatically given by (the proof of) Corollary \ref{cohcomp}: we consider as before the Lindenbaum-Tarski algebra of formulas modulo $\mathbb{T}$ and a $\kappa$-complete, $\kappa$-prime filter $\mathcal{F}$ preserving all $\kappa^+$-joins and not containing $\phi$. Then we extend the language by adding a fresh constant $c_x$ for each variable $x$ and simply set $\mathbb{T}'=\{\psi(\mathbf{c}_{\mathbf{x}}/\mathbf{x}): \psi(\mathbf{x}) \in \mathcal{F}\}$. For the second part of the proof, we consider a countable set $C_{\omega}=\cup_{i \in \omega}C_i$ of disjoint sets $C_i$ of $\kappa$-many fresh constants, and build a tree of height $\omega$ in the following way: first we start with the bottom node consisting of the constants $C_0$ in $\mathbb{T}_0:=\mathbb{T}'$, and at each level $n$ add $\kappa$-many new successor nodes of underlying sets $\cup_{i \leq n}C_i$ for each node in level $n-1$. For each pair of sentences $\phi, \psi$ over the language extended with constants from $\cup_{i \leq n}C_i$ such that $\mathbb{T}_{n-1, i}, \phi \nvdash \psi$, we consider the corresponding prime theory $\mathbb{T}_{n, j}$ as above, and put a structure in the underlying set of the node by considering the atomic formulas in $\mathbb{T}_{n, j}$. Then it is easy to see, with the same argument as in \cite{dvd}, that in the resulting Kripke model we have $(n, j) \Vdash \psi$ if and only if $\mathbb{T}_{n, j} \vdash \psi$. This construction yields, for a theory that does not prove the sentence $\phi$, a Kripke countermodel.
\end{proof}

\begin{cor}
Infinitary first-order intuitionistic logic over $\mathcal{L}_{\kappa^+, \kappa}$ with at least one constant symbol has the disjunction and existence properties. That is, if the sequents $\bigvee_{i<\kappa}\phi_i$ (respectively $\exists \mathbf{x} \phi(\mathbf{x})$) are provable, there is $i<\kappa$ and closed terms $\mathbf{t}$ such that the sequents $\phi_i$ (respectively $\phi(\mathbf{t})$) are already provable.
\end{cor}

\begin{proof}
This uses the same semantic proofs as in the case of $\mathcal{L}_{\kappa, \kappa}$, cf. \cite{espindola}.
\end{proof}

\begin{rmk}
Theorem \ref{bt} is a direct generalization of the completeness theorem for intuitionistic first-order logic over $\mathcal{L}_{\omega_1, \omega}$ of \cite{nadel}. The transfinite transitivity property for $\kappa=\omega$ is provable from the rest of the axioms (cf. \cite{espindola}), and dually, the dual distributivity rule reduces, when $\kappa=\omega$, to the law $\bigwedge_{i<\gamma} (\phi \vee \psi_i) \to \phi \vee \bigwedge_{i<\gamma} \psi_i$. Thus, Nadel's completeness theorem extends to every weakly compact cardinal.
\end{rmk}

\subsection{Completeness of first-order intuitionistic logic over $\mathcal{L}_{\kappa^+, \kappa, \kappa}$}

When restricting ourselves to a language where only conjunctions of less than $\kappa$ many formulas are allowed, Proposition \ref{filter} is enough to derive Corollary \ref{cohcomp} without needing the hypothesis of weak compactness, requiring instead only that $\kappa$ be regular and satisfies $\kappa^{<\kappa}=\kappa$. Discarding the dual distributivity rule (which involves conjunctions of $\kappa$ many formulas), we thus obtain a result that was already stated in \cite{espindola}:

\begin{thm}\label{cohcomp2}
The $(\kappa^+, \kappa, \kappa)$-coherent fragment is complete with respect to set-valued models.
\end{thm}

Finally, using Theorem \ref{cohcomp2}, we can adapt the proof of completeness for infinitary intuitionistic logic to this restricted language to get the following:

\begin{thm}\label{iilc}
 If $\kappa$ is regular and $\kappa^{<\kappa}=\kappa$, intuitionistic first-order theories over $\mathcal{L}_{\kappa^+, \kappa, \kappa}$ of cardinality at most $\kappa$ are complete with respect to Kripke models.
\end{thm}

As a final comment, it should be clear now that Theorem \ref{iilc} also sheds light on a problem left open in \cite{espindola}:

\begin{cor}\label{hc}
 Assume $GCH$. If $\kappa$ is inaccessible, intuitionistic first-order theories of cardinality less than $\kappa$ over $\mathcal{L}_{\kappa, \kappa}$ are complete with respect to Kripke models. As a consequence, Heyting cardinals are precisely the inaccessible cardinals.
\end{cor}

\section{Acknowledgements}

This research has been supported through the grants P201/12/G028 and 19-00902S from the Grant Agency of the Czech Republic.

\bibliographystyle{amsalpha}

\renewcommand{\bibname}{References} % changes the header from Bibliography to References

\bibliography{references}

%\begin{thebibliography}{widest entry}

%\end{thebibliography}

\end{document}